\pdfoutput=1

\documentclass[runningheads]{llncs}

\usepackage{graphicx}
\usepackage[short,nocomma]{optidef}
\usepackage{amsfonts}
\usepackage[caption=false]{subfig}

\begin{document}
\title{Transfer-Expanded Graphs for \\ On-Demand Multimodal Transit Systems}
\titlerunning{Transfer-Expanded Graphs for ODMTS}
\author{Kevin Dalmeijer \and Pascal Van Hentenryck}
\authorrunning{K. Dalmeijer and P. Van Hentenryck}
\institute{Georgia Institute of Technology, Atlanta, GA 30332, USA\\
\email{dalmeijer@gatech.edu, pvh@isye.gatech.edu}}
\maketitle 

\begin{abstract}
This paper considers a generalization of the network design problem for On-Demand Multimodal Transit Systems (ODMTS).
An ODMTS consists of a selection of hubs served by high frequency buses, and passengers are connected to the hubs by on-demand shuttles which serve the first and last miles.
This paper generalizes prior work by including three additional elements that are critical in practice.
First, different frequencies are allowed throughout the network.
Second, additional modes of transit (e.g., rail) are included.
Third, a limit on the number of transfers per passenger is introduced.
Adding a constraint to limit the number of transfers has a significant negative impact on existing Benders decomposition approaches as it introduces non-convexity in the subproblem.
Instead, this paper enforces the limit through transfer-expanded graphs, i.e., layered graphs in which each layer corresponds to a certain number of transfers.
A real-world case study is presented for which the generalized ODMTS design problem is solved for the city of Atlanta.
The results demonstrate that exploiting the problem structure through transfer-expanded graphs results in significant computational improvements.
\keywords{Combinatorial optimization \and Multimodal transportation \and Benders decomposition \and Transfer-expanded graphs.}
\end{abstract}

\section{Introduction}

This paper is motivated by the design and implementation of an
On-Demand Multimodal Transit System (ODMTS) for the city of Atlanta. The share of
public transit in Atlanta (about 2--3\%) is very low compared to other
American cities (e.g., about 15\% in Boston) and Atlanta is also the
8th most congested city in the world. There is thus a strong need for
a modern transit systems that leverages the train and bus
infrastructure of the city and complements it with innovative mobility
concepts.

This paper considers the design of an ODMTS for Atlanta that combines
a network of trains and buses with on-demand multimodal shuttles that
act as feeders to/from the bus/rail network and serve local
demand. ODMTS address the first/last mile problem that plagues transit
systems, while mitigating congestion on high-density corridors and
leveraging economy of scale. ODMTS and their design challenge was
introduced in \cite{MaheoEtAl2019-BendersDecompositionDesign}, which
also presents an overview of related work. The main contribution of this 
paper is to generalizes prior work by
including three additional elements that are critical for ODMTS in
large cities such as Atlanta.  First, different frequencies are
allowed throughout the network.  Second, additional modes of transit
(e.g., rail) are included.  Third, a limit on the number of transfers
per passenger is introduced.  Adding a constraint to limit the number
of transfers has a significant negative impact on existing Benders
decomposition approaches as it introduces non-convexity in the
subproblem. Instead, this paper enforces the limit through
transfer-expanded graphs, i.e., layered graphs in which each layer
corresponds to a certain number of transfers.  A real-world case study
is presented for which the generalized ODMTS design problem is solved
for the city of Atlanta.  The results demonstrate that exploiting the
problem structure through transfer-expanded graphs results in
significant computational improvements.

\section{The Generalized ODMTS Design Problem}

This section presents the generalized ODMTS design problem that
enhances the model from
\cite{MaheoEtAl2019-BendersDecompositionDesign} along several
dimensions: The choice of bus frequencies, additional transportation modes and, most importantly, a constraint on the
number of transfers. The Benders decomposition approach in
\cite{MaheoEtAl2019-BendersDecompositionDesign} exploits a natural
decomposition of the ODMTS design problem.  The network design is
determined by the \emph{master problem}, while the routing of the
passengers for a given design is determined by the \emph{subproblem}.  A
major benefit of this decomposition is that the subproblem can be
solved for each trip independently. The same decomposition is used in this paper.

\subsection{The Master Problem For Network Design}

Consider a directed multigraph $G=(V,A)$, with vertices $V =
\{1,\hdots,n\}$ and arc set $A$. Let $F$ be the set of possible
frequencies, i.e., the total number of vehicles during the time
horizon, let $M$ be the set of possible transportation modes,
which may include shuttles, and let $K$ be the total number of arcs that
each passenger may travel.  By definition, $K$ is equal to the maximum
number of transfers plus one. In the multigraph $G$, each arc $a \in
A$ is uniquely defined by the quadruple $a = (i, j, m, f) \in V \times
V \times M \times F$, $i \neq j$.  Using arc $a$ means traveling from
$i$ to $j$ with mode $m$, which departs with frequency $f$.  For a
given arc $a\in A$, these elements are referred to as $i(a)$, $j(a)$,
$m(a)$, and $f(a)$, respectively.

Designing a generalized ODMTS amounts to deciding which arcs $a \in A$ are made available to passengers.
Let the binary variable $z_a \in \mathbb{B}$ be equal to one if arc $a$ is made available, and zero otherwise.
The cost of enabling arc $a$ is given by the parameter $\beta_a$.
It is assumed that $\beta_a \ge 0$ for all $a \in A$.

For a given design, a cost is incurred due to passengers traveling trough the network.
This cost $\Phi(z)$ is a function of the values of the $z$-variables that define the design.
The value of $\Phi(z)$ can be found by solving the subproblem, which is discussed in Section~\ref{sec:subproblem}.
If the subproblem is not feasible, then $\Phi(z) = \infty$.

A formulation for the master problem is presented in Figure~\ref{fig:formulation}.
For convenience, $\delta^+(i)$ is defined as the set of all arcs going out of $i \in V$.
Similarly, the set $\delta^+(i,m)$ is defined as the set of all arcs with mode $m\in M$ going out of $i \in V$.
The sets $\delta^-(i)$ and $\delta^-(i,m)$ are defined analogously for the incoming arcs.

\newcounter{masterFormulation}
\setcounter{masterFormulation}{\value{equation}}
\addtocounter{masterFormulation}{1}
\begin{figure}[t]
	\makebox[\textwidth][c]{
		\fbox{\begin{minipage}{\textwidth}
				\begin{alignat}{3}
				~ & \min \quad \mathrlap{\sum_{a \in A} \beta_a z_a + \Phi(z),} \label{eq:objectiveOveral} \tag{\arabic{masterFormulation}a} &&&\\
				~ & \textrm{s.t.} & \sum_{a \in \delta^+(i,m)} f(a) z_a - \sum_{a \in \delta^-(i,m)} f(a) z_a & = 0 & \quad & \forall i \in V, m \in M, \label{eq:masterFlowConservation} \tag{\arabic{masterFormulation}b}\\
				~ & ~ & \sum_{f \in F \vert (i,j,m,f) \in A} z_{(i,j,m,f)} & \le 1 && \forall i \in V, j \in V, m \in M, \label{eq:masterOneFrequency} \tag{\arabic{masterFormulation}c}\\
				~ & ~ & z_a & \in \mathbb{B} && \forall a \in A. \label{eq:masterVariables} \tag{\arabic{masterFormulation}d}
				\end{alignat}
	\end{minipage}}}
	\caption{Formulation for the generalized ODMTS design problem.}
	\label{fig:formulation}
	\vspace{-0.5cm}
\end{figure}
\addtocounter{equation}{1}

Objective~\eqref{eq:objectiveOveral} minimizes the cost of the design
plus the cost of routing the passengers through the network.
Constraints~\eqref{eq:masterFlowConservation} ensure that the
frequencies for each mode are balanced at each vertex.  For example,
if three buses arrive during the time horizon, then three buses
should also depart.  Constraints~\eqref{eq:masterOneFrequency} enforce that
only one frequency can be selected for a given connection and a given
mode.  Equations~\eqref{eq:masterVariables} state the integrality
requirements.

\subsection{The Subproblem: Routing Passengers Through the Network}
\label{sec:subproblem}

For a given design, the passenger trips are routed through the network
at minimum cost.  Let $T$ be the set of all passenger trips, and let
each trip $r \in T$ be defined by an origin $o(r)$, a destination
$d(r)$, and a number of passengers $p(r)$.  If trip $r \in T$ is
routed through arc $a\in A$, then a cost of $\gamma_a^r$ is incurred.
The total cost of routing all passenger trips, $\Phi(z)$, is the sum
over the costs per trip.  It is assumed that $\gamma_a^r > 0$ for
every arc $a \in A$ and trip $r \in T$, such that the optimal route is
a simple path from $o(r)$ to $d(r)$.

Solving the subproblem amounts to solving a shortest path problem from
$o(r)$ to $d(r)$ for each trip $r \in T$, with the additional
restriction that the number of arcs in the path is at most $K$.  This
problem is known as the cardinality-constrained shortest path problem
(CSP) \cite{DahlRealfsen2000-CardinalityConstrainedShortest}.  Note
that the cardinality constraint follows from the limit on the number
of transfers.  Without this limit, the subproblem is an
(unconstrained) shortest path problem (SP), as is the case in
\cite{MaheoEtAl2019-BendersDecompositionDesign}.

It is well-known that SP possesses total unimodularity and can
be solved by linear programming (LP).  Adding an additional
constraint, however, typically destroys this structure
\cite{AnejaNair1978-ConstrainedShortestPath}.
This is indeed the case when a cardinality constraint is added to the
subproblem formulation in
\cite{MaheoEtAl2019-BendersDecompositionDesign}.  As a result, the
cost function $\Phi(z)$ would change from convex to non-convex, which
negatively impacts Benders decomposition approaches (see
Section~\ref{sec:Benders}).

{\em To remedy this limitation, this paper presents a new formulation for
the subproblem that enforces the transfer limit without destroying
total unimodularity.}  This formulation uses transfer-expanded graphs,
i.e., layered graphs with a each layer for each number of transfers. Transfer-expanded graphs encode the transit constraints 
directly, making it possible to use shortest-path algorithms.

\subsection{Transfer-Expanded Graphs}
\label{sec:transfergraphs}

Transfer-expanded graphs shares some similarities with time-expanded
networks, where each vertex has multiple copies for different periods
of time. This is the case, for example, for modern algorithms for evacuation
planning and scheduling
\cite{HasanHentenryck2017-ColumnGenerationAlgorithm,PillacEtAl2015-ColumnGenerationApproach,PillacEtAl2016-ConflictBasedPath}. Reference
\cite{MocciaEtAl2009-ColumnGenerationHeuristic} also uses a layered
network to solve the dynamic generalized assignment problem.  As a
result, some of the side-constraints do not need to be handled
explicitly. See \cite{BolandEtAl2017-ContinuousTimeService} for a
recent literature review on time-expanded graphs.

Let $\bar{G}^r = (\bar{V}^r, \bar{A}^r)$ be the
transfer-expanded graph for a given trip $r\in T$.  This graph
contains multiple copies of the original arcs and vertices, organized
in $K+1$ layers.  It is assumed that $K \ge 2$, as the subproblem is
trivial for $K=1$.  A vertex $\bar{v} = (i,k) \in \bar{V}^r$ in the
transfer-expanded graph is defined by a vertex $i \in V$ in the
original graph and by a layer $k \in \{1,\hdots,K+1\}$.  Similarly,
the definition of an arc is extended to $\bar{a} = (a, k, l)$, in
which $a \in A$ is the original arc, $k \in \{1,\hdots,K\}$ is the
layer of the starting vertex of $\bar{a}$ and $l \in \{2,\hdots,K+1\}$
is the layer of the ending vertex.

\begin{figure}[!t]
	\centering
	\includegraphics[scale=0.8]{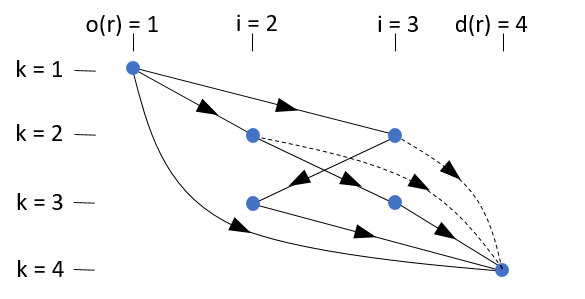}
	\caption{A transfer-expanded graph for $K=3$, $\lvert M
          \rvert = 1$, $\lvert F \rvert = 1$, for a complete graph as
          the original graph. The dotted arcs are removed when the
          triangle inequality holds.}
	\label{fig:exampleexpanded}
\end{figure}

The transfer-expanded graph is constructed as follows.  For
convenience, Figure~\ref{fig:exampleexpanded} provides an example for
$K=3$.  First, the vertex set $\bar{V}^r$ is defined.  For the origin
and the destination of the trip, introduce the vertices $(o(r),1)$ and
$(d(r),K+1)$.  For the other vertices $i \in V\backslash\{o(r),d(r)\}$
of the original graph, add the copies $(i,k)$ for $k \in \{2, \hdots,
K\}$ to the transfer-expanded graph.  The arc set $\bar{A}^r$ is
constructed based on the arcs of the original graph, as follows:

\begin{enumerate}
\item For each arc starting in the origin, i.e., $a\in
          \delta^+(o(r))$, add the arc $(a,1,2)$ if $j(a) \neq d(r)$,
          or the arc $(a,1,K+1)$ if $j(a) = d(r)$.

\item For each arc not adjacent to the origin or the destination,
  i.e., $a \in A$ and $i(a),j(a) \notin \{o(r),d(r)\}$, add the arcs
  $(a, k, k+1)$ for all $k \in \{2, \hdots, K-1\}$.

\item For each arc ending in the destination that does not start in
  the origin, i.e., $a \in \delta^-(d(r))$, $i(a) \neq o(r)$, add the
  arcs $(a, k, K+1)$ for all $k \in \{2, \hdots, K\}$.
\end{enumerate}

\noindent
By construction, it follows that solving CSP on the original graph is
equivalent to solving SP on the transfer-expanded
graph. Figure~\ref{fig:expandedsubproblem} formulates the subproblem
as a collection of SPs on transfer-expanded graphs.  Let
$y_{\bar{a}}^r \in \mathbb{B}$ be the flow on arc $\bar{a} \in
\bar{A}^r$ of trip $r \in T$.  For convenience, define
$\bar{\delta}^+_r(\bar{v})$ to be the set of all arcs in $\bar{A}^r$
coming out of $\bar{v} \in \bar{V}^r$.  Similarly, let
$\bar{\delta}^-_r(\bar{v})$ be the set of incoming arcs.

\newcounter{subproblemFormulation}
\setcounter{subproblemFormulation}{\value{equation}}
\addtocounter{subproblemFormulation}{1}
\begin{figure}[t]
	\makebox[\textwidth][c]{
		\fbox{\begin{minipage}{\textwidth}
				\begin{alignat}{3}
				\Phi(z) = & \min \quad \mathrlap{\sum_{r \in T} \sum_{\bar{a} = (a,k,l) \in \bar{A}^r} \gamma_a^r y_{\bar{a}}^r} \label{eq:extsubObjective} \tag{\arabic{subproblemFormulation}a}\\
				~ & \textrm{s.t.} & y_{\bar{a}}^r & \le z_a & \qquad\qquad & \forall r\in T, \bar{a} = (a,k,l) \in \bar{A}^r, \label{eq:extsubArcCapacity} \tag{\arabic{subproblemFormulation}b}\\
				~ & ~ & \mathllap{\sum_{\bar{a} \in \bar{\delta}^+_r(\bar{v})} y_{\bar{a}}^r - \sum_{\bar{a} \in \bar{\delta}^-_r(\bar{v})} y_{\bar{a}}^r} & \mathrlap{= \begin{cases} 1 &\textrm{ if } \bar{v} = (o(r), 1) \\ -1 &\textrm{ if } \bar{v} = (d(r), K+1) \\ 0 &\textrm{ else } \end{cases} \quad \forall r \in T, \bar{v} \in \bar{V}^r,} \label{eq:extsubFlowConservation} \tag{\arabic{subproblemFormulation}c}\\
				~ & ~ & y_{\bar{a}}^r & \ge 0 && \forall r \in T, \bar{a} \in \bar{A}^r. \label{eq:extsubVariables} \tag{\arabic{subproblemFormulation}d}
				\end{alignat}
	\end{minipage}}}
	\caption{Formulation for the subproblem on transfer-expanded graphs.}
	\label{fig:expandedsubproblem}
	\vspace{-0.5cm}
\end{figure}
\addtocounter{equation}{1}

Objective~\eqref{eq:extsubObjective} minimizes the cost of all trips.
Constraints~\eqref{eq:extsubArcCapacity} state that passengers can
only use arcs available in the design.
Constraints~\eqref{eq:extsubFlowConservation} enforce flow
conservation, and Equations~\eqref{eq:extsubVariables} define the
variables.  Due to total unimodularity of the SPs, no integrality
conditions are required.

The main advantage of using tranfer-expanded graphs is that the limit
on the number of transfers can be enforced without destroying total
unimodularity.  A potential downside is that the number of variables
and constraints in the subproblem increases with $K$.  In public
transit, however, the number of transfers that passengers are willing
to take, and therefore the value of $K$, is typically very low.
Furthermore, a larger subproblem does not necessarily mean that the
subproblem is more difficult to solve, as algorithms may benefit from
the fact that the transfer-expanded graph is acyclic.  When the
$z$-variables are integers, for example, the acyclic subproblem for
each trip can be solved in linear time through topological sorting
\cite{CormenEtAl2009-IntroductionAlgorithms}.

Finally, it is worth pointing out that if $o(r)$ and $d(r)$ are only
served by shuttles, and shuttles satisfy the triangle inequality, then
some arcs may be removed from the transfer-expanded graph without
sacrificing optimality.  Specifically, using a shuttle on the path
$(o(r),1) \rightarrow (i, 2) \rightarrow (d(r),K+1)$ for $i\in V$ is
always dominated by using a direct shuttle from $(o(r),1)$ to
$(d(r),K+1)$.  It follows that the shuttle arcs between $(i,2)$ and
$(d(r),K+1)$ may be removed for all $i\in V$, as also indicated in
Figure~\ref{fig:exampleexpanded}.  For $K\le 3$, it then follows that
the transfer-expanded graph does not require more edges than the
original graph.

\section{Benders Decomposition}
\label{sec:Benders}

Following \cite{MaheoEtAl2019-BendersDecompositionDesign}, a Benders
decomposition approach is presented for the generalized ODMTS design
problem.  The goal is to solve the master
problem~(\arabic{masterFormulation}), which is complicated by the fact
that $\Phi(z)$ is defined implicitly.  To apply Benders decomposition,
replace $\Phi(z)$ in Objective~\eqref{eq:objectiveOveral} by a new
variable $\theta \in \mathbb{R}$, and add the constraint $\theta \ge
\Phi(z)$.  Note that this does not change the problem, as $\theta =
\Phi(z)$ in any optimal solution. In Benders decomposition, the constraint $\theta \ge \Phi(z)$ is enforced through \emph{Benders cuts}.
For subproblem~(\arabic{subproblemFormulation}), these cuts are 
\begin{equation}
\label{eq:extbenderscut}
\theta \ge \Phi(\bar{z}) + \sum_{r\in T} \sum_{a \in A} \sum_{k=1}^K \lambda_a^{rk}(\bar{z}) (z_a - \bar{z}_a),
\end{equation}
with $\lambda_a^{rk}(z)$ the dual values of Constraints~\eqref{eq:extsubArcCapacity} and $\bar{z}$ any feasible solution to the LP relaxation of the master problem \cite{Benders1962-PartitioningProceduresSolving}.
For the case study in this paper, the subproblem is always feasible.
If this assumption is not satisfied, \emph{Benders feasibility cuts}, which are similar to \eqref{eq:extbenderscut}, may also be included \cite{Benders1962-PartitioningProceduresSolving}.

The Benders decomposition approach is implemented in C++ and Gurobi
8.1.1.  The master problem is the main model, and the Benders
cuts~\eqref{eq:extbenderscut} are separated in both the MIP solution
callback (in case the $z$-variables are integer) and in the MIP node
callback (in case the $z$-variables are fractional).  The subproblem
for each trip is also solved with Gurobi, and dual simplex is used to
ensure that the basis remains feasible when the subproblem is solved
for different values of $z$.  To prevent excessive calls to the
subproblem, feasibility heuristics are disabled.  The number of cut
separation rounds in the root node is set to the maximum value to get
the best possible bound.  Finally, the $2\epsilon$-trick is used to
stabilize the master problem
\cite{FischettiEtAl2017-RedesigningBendersDecomposition}.  This
stabilization uses $\epsilon = 0.00001$ and the trivial core point
obtained by assigning $z_a = \frac{1}{4}$ to every bus arc.

Without transfer-expanded graphs, the subproblem is not totally unimodular and
$\Phi(z)$ is not convex (see Section~\ref{sec:subproblem}).  In that
case, Benders decomposition cannot be applied directly.  Instead,
$\theta \ge \Phi(z)$ may be enforced by adding \emph{combinatorial
  Benders cuts} in the MIP solution callback and Benders cuts for the LP relaxation of the
subproblem in both callbacks \cite{CodatoFischetti2006-CombinatorialBendersCuts,LaporteEtAl2002-IntegerLShaped}.
However, it is well-known that relying too much on combinatorial
Benders cuts may result in slow algorithmic progress, and many cuts
may be necessary to find the optimal solution.

\section{Atlanta as a Case Study}

The generalized ODMTS design problem was solved for the city of Atlanta.  In Atlanta, the Metropolitan Atlanta Rapid Transit Authority (MARTA) operates two modes: bus and rail. The case study added on-demand shuttles and the bus
system was redesigned accordingly. More precisely, define the three modes $M=\{S,B,R\}$ for shuttle, bus, and
rail respectively.  Shuttle arcs are introduced to connect from
origins to hubs and from hubs to destinations, as well as to serve
the local demand. The corresponding $z_a$ variables are fixed to one, as
shuttles are always available. Following \cite{MaheoEtAl2019-BendersDecompositionDesign}, the cost of using a
shuttle is a weighted sum of cost and convenience, controlled by the
parameter $\alpha \in [0,1]$.  Let $d_a$ and $t_a$ be the travel
distance and the travel time of arc $a \in A$, respectively.  The
parameter $c^S$ is the cost of using a shuttle per person per unit of
distance.  The cost of traversing arc $a \in A$ for trip $r \in T$ is
then defined as $\gamma_a^r = p(r)\left((1-\alpha) c^S d_a + \alpha
t_a\right)$.

Bus arcs are defined between the potential hub locations and between
each hub and the three nearest rail stations.  The cost of enabling
bus arc $a \in A$ is given by $\beta_a = (1-a)c^Bf(a)d_a$.  That is,
the distance is multiplied by the cost per unit distance and the
number of buses over the time horizon.  The cost of traversing a bus
arc is given by $\gamma_a^r = \alpha \left(t_a + L + \frac{H}{2f(a)}
\right)$.  Here $L$ is the fixed time required for a transfer, $H$ is
the time horizon, and $\frac{H}{2f(a)}$ is the expected waiting time
before the next bus arrives, which depends on the frequency. Rail arcs are defined between all rail stations that are connected by the same rail line. The costs of traversing an arc is defined in the same way as for the buses.
For each rail arc $a \in A$, the variable $z_a$ is fixed to one, which makes the cost of enabling an arc irrelevant.

The case study uses the following data and parameters to create a realistic instance and evaluate the computational benefit of transfer-expanded graphs. It uses passenger trip data provided by MARTA for March 16, 2018, between 6am and 10am. Connecting trips have been chained together to obtain origin and destination pairs. This resulted in 2588 unique trips, with 7167 passengers in total. There are 5563 bus stops and rail stations in total, and their locations were also provided by MARTA. Eleven hubs were selected manually on the map.
For the distances $d_a$, great-circle distances are used. To estimate travel times $t_a$, the distances are divided by a constant speed of 30 mph. The cost parameters are set to $c^S = 5$ and $c^B = 1$. The fixed transfer time is chosen to be five minutes, i.e., $L=5$ minutes, and the time horizon is set to four hours, i.e., $H=240$ minutes. To balance cost and convenience, $\alpha = 0.5$ is used. The rail frequency is assumed to be fixed to six per hour, i.e., $f(a) = 6 \times 4 = 24$, and bus frequencies are determined by the model to be either three per hour or six per hour.
At most two transfers are allowed, i.e., $K=3$.

\begin{figure}[t!]
	\centering
	\subfloat[Based on transfer-expanded graphs.\label{fig:atlanta}]{%
		\centering
		\includegraphics[width=0.49\textwidth]{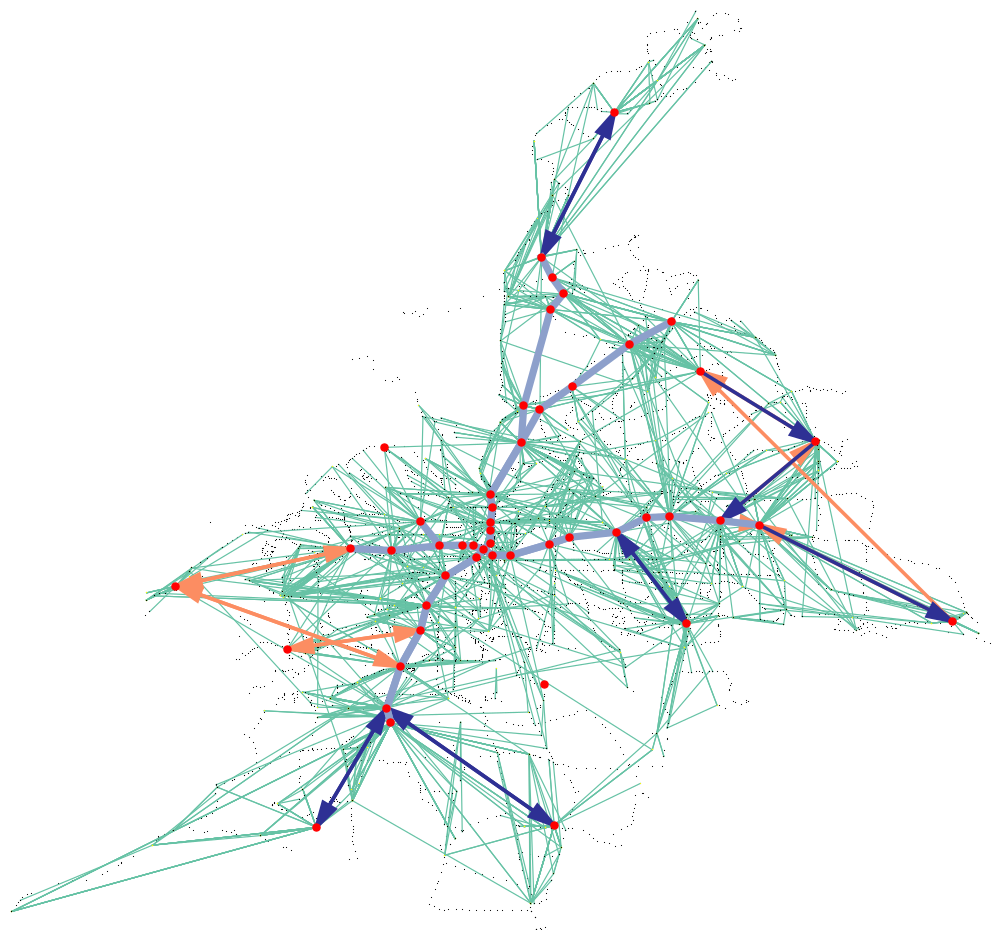}
	}
	\hfill
	\subfloat[Based on relaxed problem.\label{fig:relaxed}]{%
		\centering
		\includegraphics[width=0.49\textwidth]{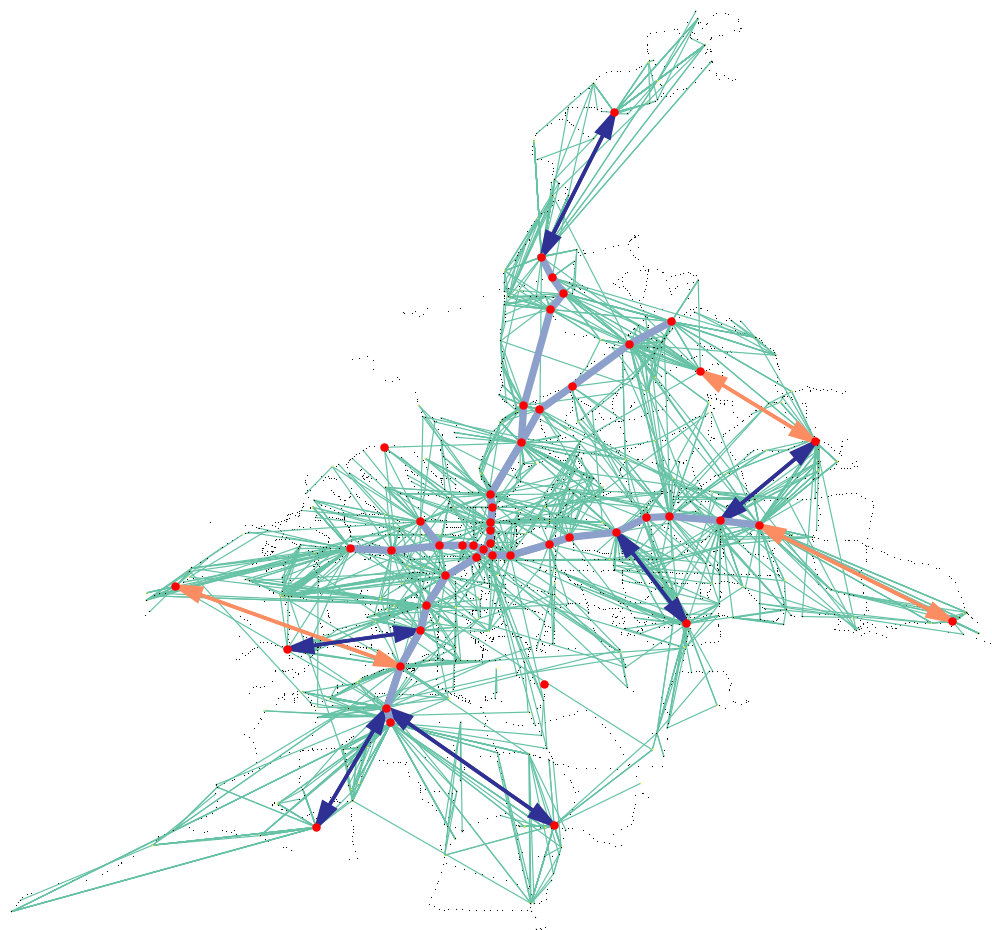}
	}
	\caption{Network designs for Atlanta showing shuttles (thin lines), rail (thick lines), and buses (arrows, orange/light for low frequency and purple/dark for high frequency).}
	\vspace{-0.5cm}
\end{figure}

Figure~\ref{fig:atlanta} presents the result of solving the generalized ODMTS problem using transfer-expanded graphs.
In total, it took 122 seconds to obtain the optimal solution and prove optimality, with a minimum objective value of 131,905. Without transfer-expanded graphs, i.e., when adding combinatorial Benders cuts, it was not possible to obtain an optimal solution in reasonable time. Instead, the evaluation considered a relaxation in which the combinatorial Benders cuts were ignored and only the Benders cuts for the LP relaxation of the subproblem were added. Solving this relaxation to optimality took 3.8 hours. Keep in mind that this relaxation explores routes that may require many transfers. To evaluate the quality of the design obtained by the relaxation, the passengers were routed through the transfer-expanded formulation with the $z$-variables fixed to their values found in the relaxation. The result is presented in Figure~\ref{fig:relaxed} and has an objective value of 131,965. Solving the relaxed problem results in a smaller public transit network because the relaxation does not completely enforce the transfer limit.

In summary, the main benefit of the transfer-expanded formulation is the significant computational benefits it provides in capturing the transfer limit. Without transfer-expanded paths, it can be optimal to fractionally select long paths that do not adhere to this constraints. These longer fractional paths likely play a role in the difference of computational performance.

\section{Conclusion}

This paper presented a generalization of the ODMTS design problem that
introduces three critical elements in practice: different
frequencies, additional transit modes, and a limit on the number of
transfers.  Transfer-expanded graphs are introduced to handle the
transfer limit without negatively impacting existing Benders
decomposition approaches.  The Atlanta case study demonstrates that
this approach is very effective, as transfer-expanded graphs
significantly improve computational performance. Exploiting the problem
structure through transfer-expanded graphs
opens the door to designing increasingly realistic networks in the
future.  One possible extension is to incorporate the capacity of the
on-demand shuttles.  As capacity of these shuttles is typically small,
expanded networks could also be used to model capacity efficiently.

\section{Acknowledgments}
This research is partly supported by NSF Leap HI proposal NSF-1854684.

\clearpage

\bibliographystyle{splncs04}
\bibliography{references}

\end{document}